\documentstyle[12pt,amsfonts,amscd]{amsart}
\newtheorem{Theorem}{Theorem}[section]
\newtheorem{Definition}[Theorem]{Definition}
\newtheorem{Proposition}[Theorem]{Proposition}

\newtheorem{Lemma}[Theorem]{Lemma}
\newtheorem{Corollary}[Theorem]{Corollary}
\theoremstyle{remark}
\newtheorem{Example}[Theorem]{Example}

\def\CC{{\Bbb C}}

\def\dim{\operatorname{dim}}

\def\be{\begin{enumerate}}
\def\ee{\end{enumerate}}
\def\bT{\begin{Theorem}}
\def\eT{\end{Theorem}}
\def\bP{\begin{Proposition}}
\def\eP{\end{Proposition}}
\def\bD{\begin{Definition}}
\def\eD{\end{Definition}}
\def\bE{\begin{Example}}
\def\eE{\end{Example}}
\def\bL{\begin{Lemma}}
\def\eL{\end{Lemma}}
\def\bC{\begin{Corollary}}
\def\eC{\end{Corollary}}

\begin{document}
\title{Holomorphic functions on subsets of ${\Bbb C}$}
\author{Buma L. Fridman and  Daowei Ma}
\begin{abstract} Let $\Gamma $ be a $C^\infty $ curve in $\Bbb{C}$ containing $0$; it becomes $\Gamma _\theta $  after rotation
by angle $\theta $ about $0$. Suppose a $C^\infty $ function $f$ can be extended
holomorphically to a neighborhood of each element of the family $\{\Gamma
_\theta \}$. We prove that under some conditions on $\Gamma $ the function $f$ is
necessarily holomorphic in a neighborhood of the origin. In case $\Gamma $ is a
straight segment the well known Bochnak-Siciak Theorem gives such a proof for \textit{real analyticity}. We also 
provide several other results related to testing holomorphy property on a
family of certain subsets of a domain in $\Bbb{C}$. \end{abstract}


\address{ buma.fridman@@wichita.edu, Department of Mathematics,
Wichita State University, Wichita, KS 67260-0033, USA}
\address{ dma@@math.wichita.edu, Department of Mathematics,
Wichita State University, Wichita, KS 67260-0033, USA}
\maketitle \setcounter{section}{-1}
\section{Introduction}

The Bochnak-Siciak Theorem [Bo,Si] states the following. Let $f\in C^\infty (D),$ $D$ is a
domain, $0\in D\subset {\Bbb R}^n.$ Suppose  $f$ is (real) analytic on every line segment
through $0$. Then $f$ is analytic in the neighborhood of $0$ (as a function
of $n$ variables). For $n=2$ this statement can be interpreted as follows. Consider the segment $
I=\{(x,y)|x\in [-1,1],y=0\},$ $I_\theta $ its rotation by angle $\theta $
about the origin. If $f$ is real analytic on each $I_\theta $ then $f$ is
real analytic in a neighborhood of the origin as a function of two
variables. Here we are interested in examining a similar statement regarding
the holomorphic property of $f$. That is if $\Gamma $ is a $C^\infty $ curve
in $\Bbb{C}$ containing $0,$ $\Gamma _\theta $ its rotation by angle $\theta 
$ about the origin, and $f$ can be extended holomorphically to a
neighborhood of each $\Gamma _\theta $, then under what condition on $\Gamma 
$ can one claim that $f$ is holomorphic in a neighborhood of $0$? For $
\Gamma $ real analytic (including $\Gamma =I$) the answer is negative, but
for some $C^\infty $ curves the answer is positive.

The questions we are examining here as well as the Bochnak-Siciak Theorem
can be considered as solving the Osgood-Hartogs-type problems; here is a quote from [ST]: ``Osgood-Hartogs-type problems ask for properties of `objects' whose restrictions to certain `test-sets' are well known''. [ST] has a number of examples of such problems. Other meaningful and interesting problems and examples of this type one can find in ([AM], [BM],  [Bo], [LM], [Ne, Ne2, Ne3], [Re], [Sa], [Si], [Zo]), and other papers. Most of the research has been devoted to consideration of formal power series and specific classes of functions of several variables as `objects' which converge (or, in case of functions, have the property of being smooth) on each curve (or subvariety of lower dimension) of a given family. The property of a series to be convergent (or, for functions, to be smooth) is then proved.

Our work in this paper is also related to another set of specific Osgood-Hartogs-type problems. The famous Hartogs theorem states that a function $f$ in ${\Bbb{C}}^n$, $n>1$, is holomorphic if it is holomorphic in each variable separately, that is, $f$ is holomorphic in ${\Bbb{C}}^n$ if it is holomorphic on every complex line parallel to an axis. So, one can test the holomorphy of a function in ${\Bbb{C}}^n$ by examining if it is holomorphic on each of the above mentioned complex lines. There is a wide area of interesting results on testing holomorphy on subsets of ${\Bbb{C}}$, specifically on curves: see [A1-A3, AG, E, G1-G3, T1, T2] and references in those articles. Some of these results assume a holomorphic extension into the inside of each closed curve in a given family, others a ``Morera-type'' property.

In this paper we also consider testing holomorphy on subsets of $\Bbb{C}$. In addition to rotations about a point (when the subset is a curve) as mentioned in the beginning, we will allow some linear transformations to be applied to these subsets. We consider a subset $S\subset \Bbb{C}$ and form a family of ``test-sets'' by considering all images of $S$ under a (small enough) subset of $\mathcal{L}$, the set of all linear holomorphic automorphisms of ${\Bbb{C}}$. We then discuss the
conditions on $S$ under which a $C^\infty $ function given in a domain will be holomorphic in that domain if
it is holomorphic on this specific family of sets. Below is a more precise
explanation. 

Let $S\subset \Bbb{C}$. We say that $f:S\rightarrow \Bbb{C}$ is
holomorphic if $f$ is a restriction on $S$ of a function holomorphic in some
open neighborhood of $S$. Let $\Bbb{L}$ be a subset of $\mathcal{L}$. 

\textbf{Definition.} The set $S$ \textit{has Hartogs property with respect to $\Bbb{L}$} (denoted $S\in~ 
\hat{H}(\Bbb L)$) \textit{if the following holds:}

\textit{Let $\Omega \subset {\Bbb{C}}$ be a domain, $f:\Omega \rightarrow 
\Bbb{C}$ a $C^\infty $ function. Suppose for any $L\in \Bbb{L},$ $f$ restricted to 
$L(S)\cap \Omega $ is holomorphic. Then $f$ is holomorphic in $\Omega $}.

The main question we are addressing here is: which sets $S$ have
Hartogs property with respect to a given set of transformations?

We will examine this question depending on $\dim (S)$ - the real Hausdorff dimension of $S
$. 

We consider three cases and provide the following answers:

1. $\dim (S)>1$. We prove that in this case $S\in \hat{H}(\Bbb T)$, where $\Bbb T$ is the group of linear translations (Theorem~\ref{S1}).

2. $\dim (S)=1$.$\,$  Such a set may or may not have Hartogs property with respect to $\Bbb T$. In addition to examples we examine explicitly the case when $S= \Gamma$ is a $C^\infty $ curve, as referred to in the beginning of this introduction. We consider the set of transformations ${\Bbb T}_1=\{\sigma\circ\tau: \sigma\in{\Bbb T}, \tau\in {\Bbb U}\}$, where 
$\Bbb{U}$ is an open subset of the group $\Bbb{C}^{*}$. Though we do not provide a complete classification of these curves we nevertheless point out the major obstacle for a curve to have Hartogs property: real analyticity. So, in this case we essentially show that if $S$ is a $C^\infty $ curve then $S\in \hat{H}({\Bbb T}_1)$ if and only
if $S$ is not analytic (for exact statements see Proposition~\ref{A1}, and Theorem~\ref{S2}).

3. $\dim (S)<1$. As in case 2 such a set may or may not have Hartogs property with respect to $\Bbb T$. We specifically examine the situation when $S$ is a sequence with one limit point (so $\dim (S)=0$), and with a reasonable restriction (a slight change of the definition of a holomorphic function on a sequence), our investigation essentially explains that $S$ has a certain Hartogs property if and only if such a sequence does not eventually end up on an analytic curve (for the precise statement see Theorem~\ref{S3} and the discussions preceding this theorem).

\section{Main Results}
{\[Case\ 1:\dim (S)>1\]$ $}

Let $S$ $\subset \Bbb{C}$. In this section we prove the following

\bT \label{S1} If $\dim (S)>1$, then $S\in \hat{H}(\Bbb T)$. \eT

The proof of this theorem follows from several statements below. For all of them $S$ is an arbitrary subset of $\Bbb C$. First we consider the following.

Let $p\in S$. A point $t$ in $T:=\{z\in\CC: |z|=1\}$ is said to be a limit direction of $S$ at $p$ if there exists a sequence $(q_j)$ in $S$ such that $\lim_j q_j=p$ and $\lim_j\tau(p,q_j)=t$, where $\tau(p,q_j):=(q_j-p)/|q_j-p|$.

\bL \label{TS}
Let $\Omega$ $\subset \Bbb{C}$ be an open set, $p\in \Omega \cap S$ and there are at least two limit directions $t_1,t_2$ of $S$ at $p$. Suppose a function $f\in C^1 (\Omega )$ is
holomorphic on $S \cap \Omega$. If $t_1\neq \pm t_2$ then $\frac{\partial f}{\partial \overline{z}}=0$ at $p$.
\eL
\begin{pf}
The derivatives of $f$ along linearly independent directions $t_1$ and $t_2$ coincide with derivatives of a holomorphic function in the neighborhood of $p$. The statement now follows from the Cauchy-Riemann equations. 
\end{pf}
\bC \label{CS} If a set $S\subset \Bbb{C}$ has a point $p$ with at least two limit
directions $t_1\neq \pm t_2$ of $S$ at $p$, then $S$ has Hartogs property with respect to $\Bbb T$.
\eC
\begin{pf}
Let $\Omega \subset \Bbb{C}$, $f$ $\in C^\infty (\Omega )$. Suppose that for any translation $L$, $f$ is holomorphic on $L(S)\cap \Omega .$ Let $z_0\in \Omega $. Pick such an $L$, that $L(p)=z_0$. Since  $f$ is holomorphic on $L(S)\cap \Omega $, and (by choice of $p$) there are at least two limit directions $t_1\neq \pm t_2$ of $L(S)\cap \Omega$ at $z_0$, then by Lemma~\ref {TS},~  $\frac{\partial f}{\partial \overline{z}}=0$ at $z_0$. So, $\frac{\partial f}{\partial 
\overline{z}}=0$ everywhere on $\Omega $, and therefore $f$ is holomorphic on $\Omega $.
\end{pf}

For a positive integer $N$ let $S_N$ be the set of points $p$ in $S$ such that $S$ has no more than $N$ distinct limit directions of $S$ at $p$. Let $M_d$ denote the Hausdorff measure of dimension $d$. Let $D(p, r)$ denote the closed disc centered at $p$ of radius $r$.

\bL \label{MS}
For $d>1$, $M_d(S_N)=0$. Hence the Hausdorff dimension of $S_N$ is $\le1$.\eL
\begin{pf} Choose a positive integer $K$ and a positive number $\epsilon$ such that 
$$ B:={2^dN\over K^{d-1}}<1,\;\; D(0,1)\cap \{q: |\tau(0,q)-1|\le\epsilon\} \subset \cup_{j=1}^KD(j/K,1/K).$$ For a positive integer $n$ let $S_N^n$ be the set of points $p$ of $S$ 
such that there exist $N$ directions $t_k$, $k=1,\dots,N$, depending on $p$, satisfying
$$D(p, 1/n)\cap S\subset \cup_{k=1}^N\{q\in\CC: |\tau(p,q)-t_k|<\epsilon\}.$$
Fix $n$ and consider a disc $D(p', r)$, where $p'\in \Bbb C$ and $r\le 1/(2n)$. If $S_N^n\cap D(p', r)$ is not empty, let $p$ be a point of this intersection. So there exist $N$ directions $t_k$, $k=1,\dots,N$, satisfying
$$D(p, 2r)\cap S_N^n\subset D(p, 2r)\cap\cup_{k=1}^N\{q\in\CC: |\tau(p,q)-t_k|<\epsilon\}.$$
The set on the right side of the above equation can be covered by $KN$ discs of radius $(2r/K)$ with centers $$p+{2rjt_k\over K},\;j=1,\dots,K,\; k=1,\dots, N.$$
Hence $D(p',r)\cap S_N^n$ can be covered by $KN$ closed discs of radius $(2r/K)$ provided $r\le 1/(2n)$. 

Now there is a positive integer $L$ such that $S_N^n$ is covered by $L$ discs of radius $1/(2n)$: $S_N^n\subset \cup_{j=1}^L D(p_j, 1/(2n))$. Each set $S_N^n\cap D(p_j, 1/(2n))$ is covered by $KN$ discs of radius $1/(nK)$. Hence $S_N^n$ is covered by $LKN$ discs of radius $1/(nK)$. For each of these smaller discs we can proceed with the similar construction. So, continuing this way we see that for any $\nu=1,2, \dots$, the set $S_N^n$ is covered by $L(KN)^\nu$ discs of radius $(1/2n)(2/K)^\nu$. It follows that $M_d(S_N^n)\le L(KN)^\nu\cdot [(1/2n)(2/K)^\nu]^d=C B^\nu$, where $C=L/(2n)^d$. Hence $M_d(S_N^n)=0$.
Since $S_N\subset \cup_{n=1}^\infty S_N^n$, we obtain $M_d(S_N)=0$.\end{pf}

{\it Proof of Theorem~\ref{S1}.} Since $\dim (S)>1$, then by Lemma~\ref{MS},  $S\setminus S_2\neq \emptyset $. Therefore there is a point $p\in
(S\setminus S_2)$, with at least two limit directions $t_1\neq \pm t_2$ of $S$ at $p$. Now by Corollary~\ref {CS}, $S$ has Hartogs property with respect to $\Bbb T$.
\qed
\newline

{\[Case\ 2:\dim (S)=1\]$ $}

The most interesting situation in this case is when $S=\Gamma$ is a curve. By using Corollary~\ref{CS} one can easily construct curves that have Hartogs property with respect to $\Bbb T$ (any broken curve (not a segment) consisting of two links and forming an angle would be such an example). On the other hand if $\Gamma$ is a real analytic curve the following statement holds.
\begin{Proposition}\label{A1} Let $\Gamma \subset \Bbb{C}$ be a real analytic curve. Then $\Gamma $
does not have Hartogs property with respect to $\mathcal{L}$.\end{Proposition} 

\begin{pf} Consider a domain $\Omega \subset \Bbb{C}$, say the unit disk, $f=\overline{
z}=x-iy$ - a nowhere holomorphic function. We prove that $f$ can be extended
holomorphically to a neighborhood of $L(\Gamma )\cap \Omega $ for any $L\in \mathcal{L}$. Without any loss of generality we may assume $L=id$, so we now consider $\Gamma \cap \Omega $. Due to the uniqueness theorem for holomorphic
functions we only need to prove the extendability of $f$ locally for any
point $z_0\in \Gamma \cap \Omega $. Again with no loss of generality we may
assume that $z_0=0$ and that near the origin $\Gamma $ is described by the equation $y=\varphi (x)$, where $\varphi (x)$ is a real analytic function. Replacing now real coordinates with $z=x+iy$ we get an
implicit equation $\frac 1{2i}(z-\overline{z})=\varphi (\frac 12(z+\overline{
z}))$, and from here one can locally recover $\overline{z}=\psi (z)$ on $\Gamma $, where $
\psi (z)$ is holomorphic near the origin. \end{pf}

We will now concentrate on smooth curves that are not analytic. We start with the following definition.  

Let $f(z)$ be a function defined on an open set $\Omega$ in the complex plane $\CC$ containing the origin. The function $f$ is said to have a Taylor series at $0$ if there is a formal power series $g(z,w)=\sum_{jk} a_{jk} z^jw^k \in \CC[[z,w]]$ such that for each nonnegative integer $n$, 
$$f(z)-\sum_{j+k\le n}a_{jk}z^j\overline z^k =o(|z|^n).$$ The Taylor series of $f$ at 0 is $g(z,\overline z)=\sum_{jk}a_{jk}z^j\overline z^k$.

We note that every $C^\infty $ function defined in the neighborhood of $0$ has a Taylor series at $0$.

Consider a curve of the form $\Gamma:=\{t+i\phi(t): 0\le t\le b\}$, where $\phi$ is a real-valued continuous function defined on the interval $[0,b]$. The function $\phi$ is said to have a Taylor series at $0$ if there exists an $h(z):=\sum_j b_jz^j\in \CC[[z]]$ such that for each nonnegative integer $n$,
$$\phi(t)-\sum_{j\le n} b_jt^j=o(|t|^n).$$
Pick an open set $\Bbb{U\subset }$ $\Bbb{C}^{*}$, and denote ${\Bbb T}_1=\{\sigma\circ\tau: \sigma\in{\Bbb T}, \tau\in {\Bbb U}\}$.
\bT \label{S2} Let $S:=\{t+i\phi(t): 0\le t\le b\}$ be a continuous curve with $\phi(0)=0$. Suppose $\phi$ has a Taylor series at 0, and for no $\lambda>0$ is $\phi$ analytic on $[0,\lambda)$. Then $S\in \hat{H}({\Bbb T}_1)$. \eT
This theorem is a corollary of Theorem~\ref{AF} below. 

First some remarks on formal power series. $\CC[[x_1, x_2, \dots, x_n]]$ denotes the set of (formal) power series 
$$g(x_1,\dots, x_n)=\sum_{k_1,\dots, k_n\ge0} a_{k_1\dots k_n} x_1^{k_1}\cdots x_n^{k_n} $$
of $n$ variables with complex coefficients. Let $g(0)=g(0,\dots,0)$ denote the coefficient $a_{0,\dots,0}$. A  power series equals 0 if all of its coefficients $a_{k_1\dots k_n}$ are equal to 0. A  power series $g\in \CC[[x_1, x_2, \dots, x_n]]$ is said to be convergent if there is a constant $C=C_g$ such that $|a_{k_1\dots k_n} |\le C^{k_1+\cdots +k_n}$ for all $(k_1,\dots,k_n)\not=(0,\dots,0)$. 
\bL \label{MT1} Let $g\in\CC [[x,y]]$ with $g'_y\not=0$, let $h\in \CC[[x]]$ be a non-zero power series with $h(0)=0$, let $E$ be a nonempty open set in the complex plane. Suppose that $g(sx, \overline sh(x))$ is convergent for each $s\in E$. Then $g$ is convergent and $h$ is convergent.
 \eL
\begin{pf} Pick $s=c\exp (i\alpha )\neq 0$, where $c=\mid s\mid $, $s\in E$. We fix $\alpha $ and since $E$
is an open set, there is a non-empty interval $(a,b)$ so for any $c\in [a,b],
$ $c\exp (i\alpha )\in E$. Replacing $x$ with $x_1\exp (-i\alpha )$ we get $
g(sx,\overline{s}h(x))=g(cx_1,ch_1(x_1))$. So, $g(cx_1,ch_1(x_1))$ converges
for all $c\in [a,b]$. Using now Theorem 1.2 from [FM] we see that $h_1(x)$ converges, implying the convergence of $h(x)$ as well. Now if $
h(x)$ is not a monomial of the form $a_1x$, we apply Theorem 1.1 from [FM]
to conclude that $g(x,y)$ is convergent as well. For the exceptional case $
h(x)=a_1x$ we need a different selection for the range of $s\in E$. Fix a
number $l>0$ and a non-zero interval $[\beta _1,\beta _2]$, such that $s=l\exp
(i\beta )\in E$ for all $\beta \in [\beta _1,\beta _2]$. Then $g(sx,
\overline{s}h(x))=g(s_1x_1,s_1^{-1}h_1(x_1))$, where $s_1=\exp (i\beta )$, $
x_1=lx$, and $h_1(x_1)={a_1}x_1$. Applying again Theorem 1.1 from [FM]
we prove the convergence of $g(x,y)$ in this case as well.
\end{pf}

\bT\label{AF} Let $f(z)$ be a continuous function defined on an open connected set $\Omega$ in the complex plane $\CC$ containing the origin, let $\Gamma:=\{t+i\phi(t): 0\le t<b\}$ be a continuous curve with $\phi(0)=0$, and let $E$ be a connected open set in the complex plane. Suppose $f$ and $\phi$ have a Taylor series at 0, that $\phi$ is analytic on $[0,\lambda)$ for no $\lambda>0$, and that for each $s\in E$ with $s\ne0$ there exists a holomorphic function $F_s$ defined in an open set $U_s$ containing $s^{-1}\Omega\cap\Gamma$ such that $f(sz)=F_s(z)$ for $z\in s^{-1}\Omega\cap\Gamma$. Then $f$ is holomorphic in the  open set $\Lambda:=\cup_{s\in E} \Gamma_s$, where $\Gamma_s$ is the connected component of $s\Gamma$ containing the origin.
\eT

\begin{pf} Let $g(z,\overline z)$ and $h(t)$ be the Taylor series at $0$ of $f$ and $\phi$ respectively. Let $\gamma(t)=t+i\phi(t)$ and $\omega(t)=t+ih(t)$. Consider an $s\in E$ with $s\ne0$. Since 
\begin{equation}\label{hyp} f(s\gamma(t))=F_s(\gamma(t))\end{equation}
for $t\in[0,b]$, we see that 
$$g(s\omega(t),\overline s(2t-\omega(t)))=F_s(\omega(t))$$
as elements in $\CC[[t]]$. Let $\psi(t)\in \CC[[t]]$ be the inverse of $\omega(t)$ so that $\omega(\psi(t))=t$. Then 
\begin{equation}\label{der}g(st,\overline s(2\psi(t)-t))=F_s(t).\end{equation}

We claim that $g_w(z,w)\equiv 0$. Suppose that is not the case. 

By Lemma~\ref{MT1}, $g(z,w)$ and $2\psi(t)-t$ are convergent. So $\psi(t)$ is convergent and $\omega(t)$ is convergent. There is a positive number $r$ such that the disk $D(0,r)\subset \Omega$, $g(z,w)$ represents a holomorphic function in $D(0,r)\times D(0,r)$, and $\psi(z)$, $\omega(z)$ represent holomorphic functions in $D(0,r)$. By (\ref{hyp}) and (\ref{der}),
\begin{equation}\label{main}
f(s\gamma(t))=g(s\gamma(t), \overline s(2\psi(\gamma(t))-\gamma(t))),\end{equation}
provided
\begin{equation}\label{region}
t\in [0,b], s\in E,  | s\gamma(t)|<r, | s(2\psi(\gamma(t))-\gamma(t))|<r.
\end{equation}

We choose an open disc $U:=D(a,v)\subset E$ with $0<v<|a|/2$, and a positive number $c<r$, such that (\ref{region}) and (\ref{main}) are satisfied for $t\in [0,c]$ and $s\in U$. Fix a $t_0\in (0,c)$. There is an $s_0\in U$ such that $g_w(s_0\gamma(t_0),w)\not\equiv 0$. Let $z_0=s_0\gamma(t_0)$. By (\ref{main}) we have, for all $t$ sufficiently close to $t_0$, that 
\begin{equation}\label{cons}f(z_0)=f({z_0\over \gamma(t)}\cdot\gamma(t))=g(z_0, \overline z_0\cdot{2\psi(\gamma(t))-\gamma(t)\over\overline{\gamma(t)}}).\end{equation}
Since $g_w(z_0,w)\not\equiv 0$, the set $\{w: g(z_0,w)=f(z_0)\}$ is discrete. Hence the function 
$$p(t):={2\psi(\gamma(t))-\gamma(t)\over\overline{\gamma(t)}}$$
is constant for all $t$ sufficiently close to $t_0$. It follows that $p(t)$ is constant on $(0, c)$. So there is a complex constant $C$ such that 
\begin{equation}\label{final}
2\psi(\gamma(t))-\gamma(t)=C \overline{\gamma(t)},\;\; 0\le t<c.
\end{equation}
Taking derivatives at 0, we obtain
$2\psi'(0)\gamma'(0)-\gamma'(0)=C\overline{\gamma'(0)}$, which forces $C=1$, since $\psi'(0)\gamma'(0)=1$ and $2-\gamma'(0)=\overline{\gamma'(0)}$.
From (\ref{final}) and $\gamma(t)+\overline{\gamma(t)}=2t$ it follows that 
\begin{equation}\label{fin}
\psi(\gamma(t))=t,\;\; 0\le t<c.
\end{equation}
The above equation implies that $\gamma(t)=\omega(t)$ for $0\le t<c$, contradicting the hypothesis that $\gamma$ is analytic on $[0,\lambda)$ for no $\lambda>0$. Therefore $g_w(z,w)\equiv 0$.

Now $g(z,w)$ does not depend on $w$, so $g\in {\Bbb C}[[z]]$, and (\ref{der}) becomes
$$g(st)=F_s(t),$$
which implies that $g$ is convergent. Hence $g$ represents a holomorphic function in $D(0,r)$ for some $r>0$. It follows from (\ref{hyp}) that 
\begin{equation}\label{mainp}
f(s\gamma(t)))=g(s\gamma(t)),\end{equation}
provided $|s\gamma(t)|<r$.
Therefore $f$ is holomorphic in the open set $Q:=D(0,r)\cap E\Gamma$.

We now prove that $f$ is holomorphic in $\Lambda$. If $0\in \Lambda$, then $0\in Q$, and we already know that $f$ is holomphic in a neighborhood of $0$. Fix a point $p\in \Lambda$, $p\not=0$. Then $p\in \Gamma_s$ for some $s\in E$, $s\not=0$, and $q:=p/s$ is a point of $\Gamma$, so $q=t_0+i\phi(t_0)$ for some $t_0\in (0,b)$. Since $\Gamma_s$ is the connected component of $s\Gamma$ containing the origin, we see that there is a $\delta\in (0, b-t_0)$ so that $\Gamma':=\{t+i\phi(t): 0\le t\le t_0+\delta\}$ satisfies that $s\Gamma'\subset \Omega$.  There is a holomorphic function $F_s$ defined in an open set $U_s\subset s^{-1}\Omega$ containing $\Gamma'$ such that $f(sz)=F_s(z)$ for $z\in \Gamma'$. Let $V_s=sU_s$ and $G_s(z)=F_s(z/s)$. Then $ s\Gamma'\subset V_s\subset \Omega$, $G_s$ is defined on $V_s$, and $G_s(z)=f(z)$ for $z\in s\Gamma'$. Choose an $\epsilon>0$ such that the disc $D:=D(s,\epsilon)$ is contained in $E$, $D$ does not contain the origin, and $D\Gamma'\subset V_s$. We now prove that $f=G_s$ in $D\Gamma'$, hence $f$ is holomorphic in $D\Gamma'$.

Consider a $u\in D$. There is a holomorphic function $G_u$ defined in an open set $V_u\subset \Omega$ containing $u\Gamma'$ such that $f(z)=F_u(z)$ for $z\in u\Gamma'$. Since $V_s\cap V_u$ contains a neighborhood of the origin, $D(0,r)\cap V_s\cap V_u$ is non-empty. By the uniqueness theorem, in the open set $D(0,r)\cap V_s\cap V_u$, the three holomorphic functions $g$, $G_s$ and $G_u$ are equal. Thus $G_s$ and $G_u$ are equal in the connected component of $V_s\cap V_u$ containing $u\Gamma'$. It follows that $f=G_s$ on $u\Gamma'$ for each $u\in D$. Thus $f=G_s$ in $D\Gamma'$, and $f$ is holomorphic in $D\Gamma'$, which is a neighborhood of $p$. Therefore $f$ is holomorphic in $\Lambda$.
\end{pf}

{\it Proof of Theorem~\ref{S2}.} Denote $\Gamma =S=\{t+i\phi(t): 0\le t\le b\}$. Let $\Omega \subset \Bbb{C}$ be a domain, $f\in C^\infty
(\Omega )$, $z_0\in \Omega $. Without any loss of generality we may  assume that $0\in \Omega $, and (since one can
use translations to move $\Gamma $) $z_0=0$. We take $E=\Bbb U$ and consider $L_s(\Gamma
)=s\Gamma $ for $s\in E$. There is a holomorphic function $G_s(z)$ in the
neighborhood of $L_s(\Gamma )\cap \Omega $ that coincides with $f$ on that
intersection. Consider $F_s(z)=G_s(sz)$. Then $f(sz)=F_s(z)$ on $
s^{-1}\Omega \cap \Gamma $. By Theorem~\ref{AF},~ $\frac{\partial f}{\partial \overline{z}}=0$ at $z_0$. So, $\frac{\partial f}{\partial 
\overline{z}}=0$ everywhere on $\Omega $, and therefore $f$ is holomorphic on $\Omega $.
\qed

{\[Case\ 3:\dim (S)<1\]$ $}
In this case an interesting situation to examine is when $S$ is a bounded sequence $(z_n)$ (and therefore $\dim (S)=0$). By using Corollary~\ref{CS} one can easily construct sequences with one limit point that have Hartogs property with respect to $\Bbb T$. On the other hand if one takes a sequence that is located on an analytic curve, and has a limit point on that curve, such a sequence will not have a Hartogs property even with respect to the entire group $\mathcal{L}$. So, a natural hypothesis here is that in order for $(z_n)$ to have Hartogs property with respect to $\mathcal{L}$ there must be no analytic curve $\Gamma$ that ${z_n}\in \Gamma$ for large $n$. However this is not true, and one can construct a counterexample. A similar statement we prove below holds, but it requires a change in the definition of a holomorphic function on a sequence. 

We will say that a function $f$ on $(z_n)$ is holomorphic if it can be extended as a holomorphic function to a {\it connected} open neighborhood of~$(z_n)$.

If a set $S=(z_n)$ has Hartogs property with respect to $\Bbb L$ and with the above definition of a holomorphic extension, we will denote that by $S\in \hat{H_0}({\Bbb L})$.  

We need another definition for the theorem below.  

Consider a sequence $(z_n)$ of complex numbers. Write $z_n=t_n+iu_n$. We assume that $t_n>0$ and $\lim z_n=0$. The sequence $(z_n)$ is said to have a Taylor series at 0 if there is an $h(z)=\sum_j b_j z^j\in \CC[[z]]$ such that 
$$u_n-\sum_{j\le k}b_j t_n^j=o(t_n^k),\;\;\; n\rightarrow \infty,$$
for each nonnegative integer $k$. Note that $h$ has real coefficients and $b_0=0$. We say that $(z_n)$ eventually lies on an analytic curve if there exists a curve $\Gamma =\{(x,y):y=\varphi (x)\}$, with $\varphi $ - real analytic
function and $\exists N$ such that $z_n\in \Gamma $ for $n\geq N$. 
\bT \label{S3} Let $S=(z_n)$, $z_1=0$, and $(z_n)$ has a Taylor series at 0 of the form $z_n\sim t_n+i h(t_n)$, where $t_n$ are positive real numbers, and $h\in \CC [[t]]$ has real coefficients. Suppose that $(z_n)$ does not eventually lie on any analytic curve. Then $S\in \hat{H_0}({\Bbb T}_2)$, where ${\Bbb T}_2=\{\sigma\circ\tau: \sigma\in{\Bbb T}, \tau\in C^*\}$.\eT
This theorem is a corollary of the following  

\bT\label{seq} Let $f(z)$ be a continuous function defined on the unit disc $D(0,1)$ in $\CC$ that has a Taylor series at 0 and let $(z_n)$ be a sequence with $z_1=0$ that has a Taylor series at 0 of the form $z_n\sim t_n+i h(t_n)$, where $t_n$ are positive real numbers, and $h\in \CC [[t]]$ has real coefficients. Suppose that $(z_n)$ does not eventually lie on an analytic curve, and that for each $s\in\CC$ with $s\ne0$ there is a holomorphic function $F_s(z)$ defined on a connected neighborhood $U_s$ of the set $Q_s:=s^{-1}D(0,1)\cap\{z_n\}$ such that $f(sz)=F_s(z)$ for $z\in Q_s$. Then $f$ is holomorphic in a neighborhood of the origin.\eT

\begin{pf} Let $g(z,\overline z)$ be the Taylor series of $f$ at 0. Let $\omega(t)=t+ih(t)$. Then 
\begin{equation}\label{hyp1}g(s\omega(t),\overline s(2t-\omega(t)))=F_s(\omega(t))\end{equation}
as elements in $\CC[[t]]$. 
Let $\psi(t)\in \CC[[t]]$ be the inverse of $\omega(t)$. 

We claim that $g_w(z,w)\equiv 0$. Suppose that is not the case. 

Similar to the proof of Theorem~\ref{AF}, we see that $h(t)$, $\omega(t)$, $\psi(t)$ are convergent, and 
\begin{equation}\label{der1}g(st,\overline s(2\psi(t)-t))=F_s(t).\end{equation} 
There is a positive number $r$ such that $D(0,r)\subset \Omega$, $g(z,w)$ represents a holomorphic function in $D(0,r)\times D(0,r)$, and $\psi(z)$, $\omega(z)$ represent holomorphic functions in $D(0,r)$.
It follows that
\begin{equation}\label{main1}
f(sz_n)=g(sz_n, \overline s(2\psi(z_n)-z_n)),\end{equation}
provided
\begin{equation}\label{region1}
|z_n|<r, | sz_n |<r, | s(2\psi(z_n)-z_n|<r.
\end{equation}

Fix $z_0\in D(0,r)$ with $z_0\ne0$ such that $g_w(z_0,w)\not\equiv0$. Then the set $\{w: g(z_0,w)=f(z_0)\}$ is discrete. 
Equation (\ref{main1}) implies that 
\begin{equation}\label{cons1}f(z_0)=f({z_0\over z_n}\cdot z_n)=g(z_0, \overline z_0\cdot
{2\psi(z_n)-z_n\over\overline z_n})=g(z_0, w_n),\end{equation}
where $w_n:=\overline z_0(2\psi(z_n)-z_n)/\overline z_n$. Since the set $\{w: g(z_0,w)=f(z_0)\}$ is discrete, and since $\lim w_n=\overline z_0$, we see that there is a positive integer $K$ such that $w_n=\overline z_0$ for $n\ge K$. Recall that $z_n\sim t_n+i h(t_n)$.
The equation $w_n=\overline z_0$ is equivalent to $\psi(z_n)=t_n$, or $z_n=\omega(t_n)=t_n+ih(t_n)$, contradicting the hypothesis that $(z_n)$ does not eventually lie on an analytic curve. Therefore $g_w(z,w)\equiv 0$.

Now $g(z,w)$ does not depend on $w$, so $g\in {\Bbb C}[[z]]$, and (\ref{der1}) becomes
$g(st)=F_s(t)$, which clearly implies that $g$ is convergent. Hence $g$ represents a holomorphic function in $D(0,r)$ for some $r>0$. Thus  
$f(sz_n)=g(sz_n)$,
provided $|sz_n|<r$.
Therefore $f$ is holomorphic in $D(0,r)$.
\end{pf}

\end{document}